\documentclass[12pt]{amsart}
\usepackage{amsmath}
\usepackage{amsfonts,amssymb,amsthm}

\newcommand{\Emin}{\mathrm{Emin\:}}
\newcommand{\Nmin}{\mathrm{Nmin\:}}

\newtheorem{theorem}{Theorem}

\theoremstyle{remark}
\newtheorem{remark}[theorem]{Remark}

\newcommand{\vol}{{\mathrm{vol}\,}}

\begin{document}

%-------------- Author entries --------------------

\title[Spherical statistics]{A note on spheres and minima}
%Article title
%\shorttitle{DRAFT} % Shortened version for
                                             % headline title

\author{Igor Rivin}

%\address{Mathematics Department, California Institute of Technology,
%Pasadena, CA 91125}

\address{Department of Mathematics, Temple University, Philadelphia}

\curraddr{Mathematics Department, Princeton University}

\email{rivin@math.temple.edu}

\thanks{The author is supported by the NSF DMS}

\date{\today}



\begin{abstract}
We write down a one-dimensional integral formula and compute large-$n$ asymptotics for 
\[
\int_{\sum_{i=1}^n x_i^2 = 1} \min_{i=1}^n|x_i| d \sigma,
\] where $d\sigma$ is the usual rotationally invariant measure on
$\mathbb{S}^{n-1}.$  The method is general, and allows to write the mean over
the sphere of an homogeneous function (Theorem \ref{homothm}) in terms of an
expectation of a function of independent, identically distributed
Gaussians. We also write down an asymptotic formula for the minimum of a large
number of identical independent positive random variables (Theorem \ref{genmin}).
\end{abstract}

\maketitle

\section*{Introduction}
Let $\mathbb{S}^{n-1}$ be the unit sphere in $\mathbb{E}^n.$ We would like to compute the following quantity:
\[
\mbox{Emin}(n)=\dfrac{1}{\vol \mathbb{S}^{n-1}} \int_{\mathbb{S}^{n-1}}
\min_{i=1}^n |x_i| d\sigma, 
\]
where $d \sigma$ is the standard measure on $\mathbb{S}^{n-1}.$ In other
words, we want to find the expected (absolute) value of the smallest
coordinate of a unit vector in $\mathbb{E}^n.$ Direct integration seems to run
into major computational difficulties, so instead we will compute $\Nmin(n),$
  which we define as as 
\emph{the expected value of the minimum of the absolute value of $n$
  independent, identically distributed random variables with mean $0$ and
  variance $1/2.$} Before we compute $\Nmin,$ we point the connection between
$\Emin$ and $\Nmin.$ First, observe that pretty much by definition,  
\[
\Nmin(n) = c_n \int_{\mathbb{E}^n} \exp\left(-\sum_{i=1}^n x_i^2\right)
\min_{i=1}^n |x_i| d x_1 \dots d x_n, 
\]
where $c_n$ is such that 
\begin{equation}
\label{norm}
c_n \int_{\mathbb{E}^n} \exp\left(-\sum_{i=1}^n x_i^2\right) d x_1 \dots d x_n
= 1. 
\end{equation}
Now we remark that $\min_{i=1}^n |x_i|$ is a $1$-homogeneous function of the
coordinates, hence we can rewrite the integral for $\Nmin$ in polar
coordinates as follows: 
\begin{equation}
\Nmin(n) = c_n \vol \mathbb{S}^{n-1}\int_0^\infty e^{-r^2} r^n \Emin{n} d r =
c_n \Emin_n \int_0^\infty e^{-r^2} r^n dr. 
\end{equation}
Since, by the obvious substition $u = r^2,$
\[\int_0^\infty e^{-r^2} r^n dr = \frac{1}{2}\int_0^\infty e^{-u} u^{(n-1)/2}
d u = \frac{1}{2}\Gamma\left(\frac{n+1}2\right).\] 
and Eq. (\ref{norm}) can be rewritten in polar coordinates as 
\[
1=c_n \vol \mathbb{S}^{n-1} \int_0^\infty r^{n-1} dr = \dfrac{c_n \vol
  \mathbb{S}^{n-1}}2 \Gamma\left(\frac{n}2\right),\] 
we see that \[\Gamma\left(\frac{n+1}2\right)\Emin(n) =
  \Gamma\left(\frac{n}{2}\right)\Nmin(n).\] 
This implies, in particular (by Stirling's formula) that 
\[\dfrac{\Nmin(n)}{\Emin(n)} \sim \sqrt{\frac{n+1}2}.\]

It is clear that the above argument only depends on the homogeneity of the
function, so it immediately generalizes to the following:

\begin{theorem}
\label{homothm}
Let $f(x_1, \dots, x_n)$ be a homogeneous function on $\mathbb{E}^n$ of degree
$d$ (in other words, $f(\lambda x_1, \dots, \lambda x_n) = \lambda^d f(x_1,
\dots, x_n)$.) Then
\[
\dfrac{\Gamma\left(\frac{n+d}2\right)}{\vol \mathbb{S}^{n-1}}
\int_{\mathbb{S}^{n-1}} f d \sigma = 
  \Gamma\left(\frac{n}{2}\right)\mathbb{E}\left(f(\mathbf{X}_1, \dots,
  \mathbf{X}_n)\right),\]  
where $\mathbf{X}_1, \dots, \mathbf{X}_n$ are independent random variables
with probability density $e^{-x^2}.$
\end{theorem}

Now, to complete our computation, we must compute $\Nmin(n).$ Let $X_1, \dots,
X_n$ be independent, identically distributed variables, whose common
distribution $F$ is supported on $[0, \infty).$ What is the distribution of
  $X_{(1)} = \min(X_1, \dots, X_n)?.$ 
The probability that $X_{(1)}$ is greater than $y$ is obviously $(1-F(y))^n,$
so the distribution function of $X_{(1)}$ is obviously $1-(1-F(y))^n.$ It
follows (by integration by parts) that the expectation of $X_{(1)}$ is  
\[E_{(1)} = \int_0^\infty (1-F(y))^n d y.\]
In our particular case of $X_i = |Z_i|,$ where $Z_i$ is normal with variance
$1/2,$  
\begin{equation}
\label{distrib}
F(y) =  \dfrac{1}{\sqrt{\pi}} \int_{-y}^y \exp{-x^2} d x.
\end{equation}
We thus have an integral formula for $\Emin(n)$ promised in the abstract:

\begin{equation}
\label{intform}
\Emin(n) = \dfrac{\Gamma\left(\frac{n}{2}\right)}
{\Gamma\left(\frac{n+1}{2}\right)}\int_0^\infty\left(1-\dfrac{1}{\sqrt{\pi}}
  \int_{-y}^y \exp{-x^2} d x\right)^n. 
\end{equation}
Since this formula is somewhat unwieldy, it is worthwhile to state an asymptotic result. Then, 
\begin{theorem}
Let $F$ be as in Eq. (\ref{distrib}). Then
\[\int_0^\infty (1-F(y))^n d y = \dfrac{\sqrt{\pi}}{2(n+1)}+o(n^{-2}).\]
\end{theorem}
\begin{proof}
We write
\begin{equation}
\label{splitup}
\int_0^\infty (1-F(y))^n d y = I_1 + I_2 + I_3,
\end{equation}
where
\begin{gather}
I_1 = \int_0^{n^{-3/4}} (1-F(y))^n dy,\\
I_2 = \int_{n^{-3/4}}^C (1-F(y))^n d y,\\
I_3 = \int_C^\infty (1-F(y))^n d y.
\end{gather}
First, we show that $I_1$ and $I_2$  are asymptotically negligibly
small. Indeed, 
$1-F(y) = \frac{2}{\sqrt{\pi}}\int_y^\infty \exp(-x^2) d x.$ For $x > 1,$
$\exp(-x^2) \ll \exp(-x),$ and  
so $1-F(y) \ll \exp(-x),$ so for a suitable choice of $C,$ $I_3$
decreases exponentially in $n.$ Furthermore $1-F(y)$ is monotonically
decreasing, so to estimate $I_2,$ we write 
\[\int_{n^{-3/4}}^C (1-F(y))^n \leq C (1-F(n^{-3/4}))^n \ll \exp(-C_1
n^{1/2}),\]
since for small $y,$ $F(y) \approx 2 y/\sqrt{\pi}.$ 

Finally, to estimate the first integral, we expand $F(y)$ in a Taylor series,
to obtain  
\[1-F(y) = 1-\frac{2 y}{\sqrt{\pi}} + O(y^3) = \left(1-\frac{2 y}{\sqrt{\pi}}
\right)(1+O(y^3)),\] 
so that 
\[(1-F(y)]^n = \left(1-\frac{2y}{\sqrt{\pi}}\right)^n (1+O(y^3))^n.\]
for $y < n^{-3/4},$ we know that $(1+O(y^3))^n - 1 = O(n^{-5/4}),$ so that
\[\int_0^{n^{-3/4}}(1-F(y)^n d y =
\left(1+O(n^{-5/4})\right)\int_0^{n^{-3/4}}
\left(1-\frac{2y}{\sqrt{\pi}}\right)^n dy.\] 
Now, since 
\[\int_{n^{-3/4}}^{\sqrt{\pi}/2}\left(1-\frac{2y}{\sqrt{\pi}}\right)^n
dy\ll\exp\left(-\dfrac{\sqrt{\pi}}{2} n^{1/4}\right),\] it follows that
\[
\int_0^{n^{-3/4}}\left(1-\frac{2y}{\sqrt{\pi}}\right)^n dy\sim
  \int_0^{\sqrt{\pi}/2} \left(1-\frac{2y}{\sqrt{\pi}}\right)^n dy =
  \dfrac{\sqrt{\pi}}{2(n+1)}.
\] 
\end{proof}

It is clear that in the above argument we don't actually need $F$ to be the
normal distribution, and it holds in much greater generality:
\begin{theorem}
\label{genmin}
Let $X_1, \dots, X_n$ be identically independently distributed variables on
$[0, \infty]$ with distribution function $F.$ Suppose that the distribution
$F$ satisfies the following conditions:
\begin{enumerate}
\item $F$ has a continuous \emph{novanishing} density $f$ in a neighborhood of
  $0.$ 
\item $1-F \in L^p([0, \infty)),$ for some $p>0.$
\end{enumerate}
Then, as $n\rightarrow \infty,$ 
\[\mathbb{E}(\min(X_1, \dots, X_n)) \sim \dfrac{1}{f(0)(n+1)}.\]
\end{theorem}
\begin{proof} The argument goes through pretty much as above, except for the
  proof that the integral $I_3$ decreases exponentially with $n.$ This,
  however, is easily fixed: The function $1- F(y)$ is monotonically
  decreasing, so, for $y > C,$ $(1-F(y))^n \leq C^{n-p} (1-F(y))^p,$ where $p$
  is as in the statement of the theorem. The result follows immediately.
\end{proof}
\begin{remark}
The second condition in the statement of Theorem \ref{genmin} can be
interpreted as saying that if the expectation of the minimum of 
\emph{some} number of variables is finite, then we have the claimed
asymptotics, and otherwise the expectation is always infinite, so, in a sense,
we have complete asymptotic information in that case also.
\end{remark}
\begin{remark}
Theorem \ref{genmin} has no information about the error term, but this is due
to the very weak regularity assumption on the density $f$ at the origin.
\end{remark}
\bibliographystyle{amsplain}
%\begin{thebibliography}{999999999999}
%\end{thebibliography}
\end{document}